\documentclass[12pt]{article}

\usepackage{amsmath,amsfonts,amssymb}
\usepackage{graphicx}
\usepackage{geometry}
\usepackage{lmodern}
\usepackage{caption}
\usepackage{fancyhdr}
\usepackage{lastpage}
\usepackage{bm}
\usepackage{mathrsfs} 
\usepackage{tikz}
\usepackage{float}
\usepackage{enumitem}
\usepackage{placeins}
\usepackage{pgfplots}
\usepackage{booktabs}
\usepackage{lipsum}
\usepackage{graphicx}
\usepackage{xpatch}
\usepackage[colorlinks=true,allcolors=blue]{hyperref}
\usepackage[exponent-product=\cdot]{siunitx}
\usepackage{tikz}
\usetikzlibrary{arrows.meta}
\usepackage{mathtools}
\usepackage{subcaption}
\DeclareMathOperator{\mydiv}{div}
\newcommand{\A}{\alpha}
\title{\textbf{\LARGE An Air-Gap Element for the Isogeometric 
Space-Time-Simulation of Electric Machines}}
\author{\normalsize M.~Reichelt$^{1}$, M.~Wiesheu$^{2}$, M.~Merkel$^{2}$, S.~Schöps$^{2}$, and O.~Steinbach$^{1}$ \\ 
\small $^{1}$Institute of Applied Mathematics, \\
\small TU Graz, Steyrergasse 30, 8010 Graz, Austria \\
\small $^{2}$Institute for Accelerator Science and
Electromagnetic Fields (TEMF),\\
\small TU Darmstadt, Schloßgartenstr. 8, 64289 Darmstadt, Germany \\
\small E-mail: michael.reichelt@tugraz.at}
\date{}

\renewenvironment{abstract}
 {\par\noindent\textbf{\abstractname \\}\ \ignorespaces}
 {\par\medskip}

\begin{document}

\maketitle
\thispagestyle{empty}
\begin{abstract}
Space-time methods promise more efficient time-domain simulations, in particular of electrical machines. However, most approaches require the motion to be known in advance so that it can be included in the space-time mesh. To overcome this problem, this paper proposes to use the well-known air-gap element for the rotor-stator coupling of an isogeometric machine model. 
First, we derive the solution in the air-gap region and then employ it to couple the rotor and stator. This coupling is angle dependent and we show how to efficiently update the coupling matrices to a different angle, avoiding expensive quadrature. Finally, the resulting time-dependent problem is solved in a space-time setting.
The spatial discretization using isogeometric analysis is particularly suitable for coupling via the air-gap element, as NURBS can exactly represent the geometry of the air-gap. Furthermore, the model including the air-gap element can be seamlessly transferred to the space-time setting.
However, the air-gap element is well known in the literature. The originality of this work is the application to isogeometric analysis and space-time.

\end{abstract}

\section{Introduction}
In the design of electric machines, simulation plays an important role in predicting machine performance and guiding design choices. Conventionally, these simulations are performed using the Finite Element Method (FEM) with low-order elements \cite{Monk_2003aa}. When considering rotating electric machines,
space-time finite element methods \cite{Gangl2024} can resolve the air-gap between the rotor and the stator by using simplicial
space-time finite elements. However, an accurate representation
of the cylindrical geometry  poses challenges as such shapes are costly, e.g., in terms of number of elements, in particular
in the space-time approach. Since rotor and stator are well
separated, and not changing in time, it is sufficient
to consider fixed spatial meshes for both, which are coupled via the
air-gap domain within a time stepping approach. To avoid remeshing, techniques such as locked step, sliding surface,  moving band and mortar methods have been developed which introduce a (non-conforming) interface to couple rotating subdomains \cite{Rodger_1990aa,Davat_1985aa,De-Gersem_2004af}. In particular, the mortar setting was recently extended to be applied in the isogeometric setting, \cite{Bontinck_2018ac,Kapidani_2022aa}.
Another approach to realize rotation is to replace the air-gap region by another field description, e.g. the Boundary Element Method (BEM) \cite{Kurz_1999aa} or an analytical formulation, i.e., the air-gap method \cite{Abdel-Razek_1982aa}. 

Isogeometric analysis (IGA) is a well suited FEM variant for simulating electric machines, particularly due to its ability to represent geometries exactly. Introduced by Hughes et al. \cite{Hughes_2005aa}, IGA bridges the gap between computer-aided design (CAD) and finite element analysis by using CAD basis functions, i.e. B-splines and non-uniform rational B-splines (NURBS), for both geometry representation and the basis functions for the simulation. It allows for an exact geometric representation of cylindrical structures and the smoothness of IGA basis functions can improve the accuracy of simulations for the (smooth) magnetic fields in electric machine analysis.

The combination of the air-gap method with IGA is especially well suited as both methods represent the rotor and stator interfaces of the air-gap geometrically exactly.
This paper introduces, for the first time, a method that integrates a flexible rotor movement enabled by a scaled air-gap element in the IGA context, combining those approaches to overcome these difficulties and enhance simulation efficiency.
However, a naive implementation of the air-gap element leads to badly conditioned algebraic systems. Therefore, we introduce a new scaling which is inspired by \cite{Gangl_2020ag}.

The paper is structured as follows. Sections \ref{sec:em} and \ref{sec:airgap} introduce the model and in particular the air-gap element. Sections~\ref{eq:discrete} and \ref{sec:time} discuss the discretization in space and time. The quantity of interest (torque) is defined in Section~\ref{sec:torque}. A machine model is introduced in Section~\ref{sec:machine} and the simulation results are shown in Section~\ref{sec:results}. The paper closes with conclusions in Section~\ref{sec:conclusions}.

\section{Electromagnetic Model}\label{sec:em}
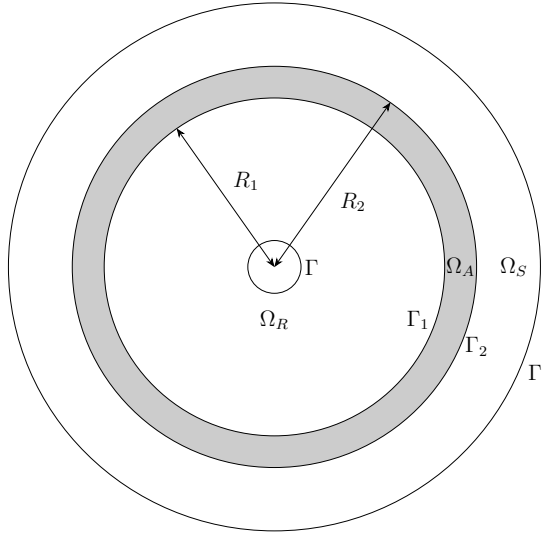
\begin{figure}[]
  \begin{center}
  \scalebox{0.7}{%
  \begin{tikzpicture}[scale=1]
    \draw (0,0) circle (5cm);
    \node at (4.5,0) {$\Omega_S$};

    \draw[fill=black,fill opacity=0.2] (0,0) circle (3.8cm);
    \node at (3.5,0) {$\Omega_A$};

    \draw[fill=white] (0,0) circle (3.2cm);
    \node at (0,-1) {$\Omega_R$};

    \draw (0,0) circle (0.5cm);
    \node at (0.7,0) {$\Gamma$};

    \draw[Stealth-Stealth] (0,0) -- (125:3.2) node [pos = 0.5, anchor = south west] {$R_1$};
    \draw[Stealth-Stealth] (0,0) -- (55:3.8) node [pos = 0.5, anchor = north west] {$R_2$};

    \node at (2.7, -1) {$\Gamma_1$};
    \node at (3.8, -1.5) {$\Gamma_2$};
    \node at (4.9, -2) {$\Gamma$};
\end{tikzpicture}
  }
  \caption{Schematic cross-section of an electric machine. {Here, the inner domain $\Omega_R$ is the rotor, $\Omega_S$ denotes the region of the stator and $\Omega_A$ denotes the region of the air-gap.}}\label{fig:geom}
  \end{center}
\end{figure}

Two-dimensional simulations of electric machines, e.g.~\cite{Salon_1995aa}, deal with complex variants of the principal geometry $\Omega$ depicted in Figure \ref{fig:geom}, where the interior domain $\Omega_R$ is the rotor (e.g. including permanent magnets), $\Omega_S$ denotes the region of the stator (e.g. containing slots with the windings), and $\Omega_A$ denotes the region of the air-gap. We employ the isogeometric variant of the finite element method (IGA), \cite{Hughes_2005aa}, i.e., those domains are exactly represented by patches of spline-based mappings. 

An appropriate low-frequency model of Maxwell's equations \cite{Jackson_1998aa} is the eddy current approximation \cite{Schmidt_2008aa}, which reads in two dimensions
\begin{align}
  \sigma \frac{\mathrm{d}}{\mathrm{d}t} A_z - \mydiv(\nu \nabla A_z) &= J_z \quad \text{in } \Omega, \\
  A_z &=0 \quad \text{ on } \Gamma,
\end{align}
with conductivity $\sigma$, reluctivity $\nu$ and the $z$-components of  magnetic vector potential $\textbf{A}$, excitation current density $\mathbf{J}$ and suitable boundary conditions. As there is no excitation current density in the air-gap, we split the domain $\Omega$ into non-overlapping subdomains, yielding
\begin{align}
    \label{eq:fem}
    \sigma \frac{\mathrm{d}}{\mathrm{d}t} A_z -\mydiv(\nu \nabla A_z) &= J_{z} \quad &\text{in } \Omega_R \cup \Omega_S,\\
    -\mydiv(\nu \nabla A_z) &= 0 \quad &\text{in } \Omega_A,
    \label{eq:ana}
\end{align}
with the constraint that $A_z$ and its conormal derivative are continuous across the interface boundaries $\Gamma_1$ and $\Gamma_2$. If we denote the solutions in the rotor, air and stator region by $A_z^R, A_z^A, A_z^S$, respectively, the transmission conditions are
\begin{align}
    A_z^R = A_z^A, \quad \nu_R \partial_{n_R} A_z^R + \nu_A \partial_{n_A} A_z^A = 0 \quad \text{on } \Gamma_1, \\
    A_z^S = A_z^A, \quad \nu_S \partial_{n_S} A_z^S + \nu_A \partial_{n_A} A_z^A = 0 \quad \text{on } \Gamma_2,
\end{align}
where the subscript of $\nu$ indicates the possibly different material behavior of the respective regions towards the interface.
If the rotor moves, the transmission conditions are angle dependent. Note that the normal vectors $n_R$, $n_S$ and
$n_A$ are defined exterior to the respective domains.

\section{Air-Gap Solution}\label{sec:airgap}
For the solution in the air-gap region we follow the ideas of \cite{Abdel-Razek_1982aa}. Due to the simple geometry, the air-gap solution of \eqref{eq:ana} can be computed analytically. 
In polar coordinates we have
\begin{align}
  A_z|_{\Omega_A}(r,\varphi) &= \A_0 + \A_0' \ln(r) \nonumber \\
  &+ \sum_{k=1}^\infty \cos(k \varphi) \left( \alpha_k r^k + \alpha_k' r^{-k} \right)+ \sum_{k=1}^\infty \sin(k \varphi) \left( \beta_k r^k + \beta_k' r^{-k} \right), \label{eq:ana_sol} 
\end{align}
where the coefficients $\A_0, \A_0', \alpha_k, \alpha_k', \beta_k, \beta_k'$ are to be determined. As there is no time derivative in \eqref{eq:ana}, and hence the solution is only dependent on boundary conditions at every point in time, for better readability we suppress the time dependency of the coefficients in this section. We impose Dirichlet boundary conditions 
\begin{align}
    A_z|_{\Omega_A} &= v_i \quad \text{ on } \Gamma_i, \quad \text{ for } i\in\{1,2\}
\end{align}%
on both interfaces, where $v_i$ is assumed to be given. For comparison of coefficients the boundary conditions are evolved into a Fourier series, i.e.
\begin{align}
    v_i = v_{0,i} + \sum\limits_{k=0}^\infty \Big[ v_{k,i} \cos\left(k \varphi \right) + v_{k,i}' \sin\left(k \varphi \right) \Big],
\end{align}
with
\begin{align}
    v_{0,i} &= \frac{1}{2 \pi R_i} \int_0^{2\pi} v_i(R_i \varphi) \, d\varphi, \label{eq:fouriercoef1}\\
    v_{k,i} &= \frac{1}{\pi R_i} \int_0^{2\pi} \cos\left( k \varphi \right) v_i(R_i \varphi) \, d \varphi,  \label{eq:fouriercoef2} \\
    v_{k,i}' &= \frac{1}{\pi R_i} \int_0^{2\pi} \sin\left(  k \varphi \right) v_i(R_i \varphi) \, d \varphi,  \label{eq:fouriercoef3}
\end{align}
where $R_i$ denotes the radius of the interface and for brevity of notation we make the values of the boundary data dependent on the arc length on the respective interface. With the Fourier coefficients \eqref{eq:fouriercoef1}--\eqref{eq:fouriercoef3} of the boundary data, enforcing Dirichlet conditions for the zero-th mode leads to the $2\times 2$ system
\begin{align}
    G_0\underline{\alpha_0}{\coloneq}
    \begin{pmatrix}
        1 & \ln\left( R_1 \right) \\
        1 & \ln\left( R_2 \right)
    \end{pmatrix}
    \begin{pmatrix}
        \A_0 \\
        \A_0'
    \end{pmatrix}
    = 
    \begin{pmatrix}
        v_{0,1} \\
        v_{0,2}
    \end{pmatrix} \label{eq:TwoByTwoConst}
\end{align}
and for all higher modes, we have to solve the two independent $2\times 2$ systems
\begin{align}
    G_k\underline{\alpha_k}{\coloneq}
    \begin{pmatrix}
        R_1^k & R_1^{-k} \\
        R_2^k & R_2^{-k}
    \end{pmatrix}
    \begin{pmatrix}
        \alpha_k \\
        \alpha_k'
    \end{pmatrix}
    =
    \begin{pmatrix}
        v_{k,1} \\
        v_{k,2}
    \end{pmatrix} \text{} \label{eq:TwoByTwoCos}
\end{align}
and
\begin{align}
    G_k\underline{\beta_k}{\coloneq}
    \begin{pmatrix}
        R_1^k & R_1^{-k} \\
        R_2^k & R_2^{-k}
    \end{pmatrix}
    \begin{pmatrix}
        \beta_k \\
        \beta_k'
    \end{pmatrix}
    =
    \begin{pmatrix}
        v_{k,1}' \\
        v_{k,2}'
    \end{pmatrix}. \label{eq:TwoByTwoSin}
\end{align}
If one interface is rotated by an angle $\delta$ the calculation of the boundary data's Fourier coefficients changes and one can use properties of trigonometric functions to obtain
\begin{align}
    \tilde{v}_{k,i}(\delta) &= \frac{1}{\pi R_i} \int_0^{2\pi} \cos\left( k \varphi \right) v_i\left(R_i (\varphi-\delta)\right) \, d \varphi \nonumber \\
    &= \frac{\cos(k\delta)}{\pi R_i} \int_0^{2\pi} \cos\left( k \varphi \right) v_i\left(R_i \varphi\right) \, d \varphi - \frac{\sin(k\delta)}{\pi R_i} \int_0^{2\pi} \sin\left( k \varphi \right) v_i\left(R_i \varphi\right) \, d \varphi \label{eq:rotation1}
\end{align}
and analogously
\begin{align}
    \tilde{v}_{k,i}'(\delta)
    &= \frac{\sin(k\delta)}{\pi R_i} \int_0^{2\pi} \cos\left( k \varphi \right) v_i\left(R_i \varphi\right) \, d \varphi  + \frac{\cos(k\delta)}{\pi R_i} \int_0^{2\pi} \sin\left( k \varphi \right) v_i\left(R_i \varphi\right) \, d \varphi \label{eq:rotation2}
\end{align}
from the data $v_i$ in the reference domain. So the rotated coefficients are just weighted sums of the Fourier coefficients in the reference domain. Noticing, that the integrals on the right hand side are the standard Fourier coefficents, this can be put into matrix form yielding
\begin{align}
    \begin{pmatrix}
        \tilde{v}_{k,i} \\
        \tilde{v}_{k,i}'
    \end{pmatrix} = 
    \begin{pmatrix}
        \cos(k\delta) & - \sin(k \delta) \\
        \sin(k\delta) & \cos(k \delta)
    \end{pmatrix}
    \begin{pmatrix}
        v_{k,i} \\
        v_{k,i}'
    \end{pmatrix},\label{eq:fourier_rotation}
\end{align}
where the occurring matrix is a rotation matrix for the angle $k\delta$.

\section{Spatial Discretization and Coupling}\label{eq:discrete}
We consider the remaining equation \eqref{eq:fem} on the reference domain with angle dependent transmission conditions. As in this viewpoint, the observer is fixed to the respective domain, the total time derivative becomes a partial one. Using a variational approach in space and applying integration by parts, yields the problem to find $A_z \in H^1_\Gamma(\Omega_D)$, such that
\begin{align}
    \int\limits_{\Omega_D} \left(\sigma \partial_t A_z V 
    +  \nu \nabla A_z \cdot \nabla V \right)\, dx 
    + \int\limits_{\Gamma_1 \cup \Gamma_2} \left(S (\delta) A_z\right) V \, ds = \int\limits_{\Omega_D} J_z V \, dx \label{eq:var_formulation}
\end{align}
for all $V \in H^1_\Gamma(\Omega_D)$, where $\Omega_D = \Omega_R \cup \Omega_S$ and $ H^1_\Gamma(\Omega_D)$ denotes $H^1$-functions with vanishing trace on $\Gamma$. $S$ is the rotation dependent Dirichlet-to-Neumann mapping, which assigns to every function in $ H^1_\Gamma(\Omega_D)$ and angle $\delta$ the according Neumann trace arising from the air-gap solution.  We discretize \eqref{eq:var_formulation} by a Galerkin approach, i.e.
\begin{align}
    A_{z,h}(t,x) = \sum\limits_{n=1}^N A_n(t) \psi_n(x)
\end{align}
where $\lbrace \psi_n \rbrace_{n=1}^N$ is the basis of the IGA space comprising rotor and stator, yielding
\begin{align}
    \sum\limits_{n=1}^N\int\limits_{\Omega_D} \left(\sigma \partial_t A_n \psi_n \psi_j +  A_n \nu \nabla \psi_n \cdot \nabla \psi_j \right)\, dx + 
    \sum\limits_{n=1}^N
    \int\limits_{\Gamma_1 \cup \Gamma_2} A_n \left(S(\delta) \psi_n \right) \psi_j \, ds 
    = \int\limits_{\Omega_D} J  \psi_j \, dx  
\end{align}
$\forall j=1,\ldots, N$. We collect all IGA degrees of freedom in a vector $\underline{A}$. The discretized Dirichlet-to-Neumann mapping $S_h(\delta)$ is realized via the air-gap solution. For that purpose, we abort the infinite series in \eqref{eq:ana_sol} at a finite number $K$. The necessary integrals for the Fourier coefficients are computed by
\begin{align}
    v_{i,0} &= \frac{1}{2\pi R_i} \int_0^{2\pi}  A_{z,h}(R_i \varphi) \, d \varphi = \sum\limits_{n=1}^N \frac{A_n}{2\pi R_i} \int_0^{2\pi} \cos\left( k \varphi \right) \psi_n(R_{i} \varphi) \, d \varphi, \\
    v_{i,k} &= \frac{1}{\pi R_i}\int_0^{2\pi} \cos\left( k \varphi \right) A_{z,h}(R_i \varphi) \, d \varphi = \sum\limits_{n=1}^N \frac{A_n}{\pi R_i}\int_0^{2\pi} \cos\left( k \varphi \right) \psi_n(R_{i} \varphi) \, d \varphi, \\
    v_{i,k}' &= \frac{1}{\pi R_i} \int_0^{2\pi} \sin\left( k \varphi \right) A_{z,h}(R_i \varphi) \, d \varphi = \sum\limits_{n=1}^N \frac{A_n}{\pi R_i} \int_0^{2\pi} \sin\left( k \varphi \right) \psi_n(R{_i} \varphi) \, d \varphi. \label{eq:FourierSum}
\end{align}
If we now collect all Fourier coefficients $v_{i,k}$ and $v_{i,k}'$ for both interfaces into a vector $\underline{F}$, then there is a linear relation between $\underline{F}$ and $\underline{A}$, i.e.
\begin{align}
    \underline{F} = C_2 \underline{A},
\end{align}
where the row $\underline{r}^\top$ of $C_2$ corresponding to e.g. \eqref{eq:FourierSum} is given by
\begin{align}
    \underline{r}[n] = \frac{1}{\pi R_i} \int_0^{2\pi} \sin\left( k \varphi \right) \psi_n(R_i \varphi) \, d \varphi.
\end{align}
The other rows follow {similarly}.
Note, that due to \eqref{eq:fourier_rotation}, the associated matrix can be efficiently adapted to different rotation angles and quadrature is only needed once for the reference domain, yielding the angle dependent matrix $C_2(\delta)$. The evaluation of $\partial_n A_z$ at the interfaces, taking the air-gap solution to determine the Neumann trace, yields similar integrals and a matrix denoted by $C_1(\delta)$. {Note, that these matrices have the same structure as the coupling matrices in \cite{Egger_2022ab}}. 
When assembling the degrees of freedom of the IGA solution $\underline{A}$ together with the air-gap solution coefficients $\underline{\alpha} = (\alpha_0,\alpha_0', \ldots, \beta_1, \beta_1' \ldots)^\top$ into a vector
\begin{align}
    \underline{U} = \begin{pmatrix}
        \underline{A} \\
        \underline{\alpha}
    \end{pmatrix}
\end{align}
we obtain the following block system of algebraic equations
\begin{align}
    \begin{pmatrix}
        M & 0 \\
        0 & 0 
    \end{pmatrix}
    \dot{\underline{U}} + \begin{pmatrix}
        K & C_1(\delta) \\
        C_2(\delta) & K_F
    \end{pmatrix}
    \underline{U}
    = \begin{pmatrix}
        \underline{J} \\
        0
    \end{pmatrix} \label{eq:dae_system}
\end{align}
which is a differential algebraic system in time. $M$ and $K$ are the standard mass and stiffness matrices of the IGA discretization, $\underline{J}$ is the standard load vector, $K_F$ comprises the $2\times2$ equation systems for the air-gap solution. $C_1(\delta)$ and $C_2(\delta)$ are defined above. Note, that the discrete Dirichlet-to-Neumann mapping
\begin{align}
    S_h(\delta) = C_2(\delta) K_F^{-1} C_1(\delta)
\end{align}
is implicitly contained in the DAE system.
For brevity we define 
\begin{align*}
    \tilde{M} = \begin{pmatrix}
        M & 0 \\
        0 & 0 
    \end{pmatrix},~
    \tilde{K} = \begin{pmatrix}
        K & C_1(\delta) \\
        C_2(\delta) & K_F
    \end{pmatrix}
    ,~
     \tilde{\underline{J}} = \begin{pmatrix}
        \underline{J} \\
        0
    \end{pmatrix}.
\end{align*}

\subsection{Conditioning}
\begin{figure}
    \centering
    \begin{subfigure}[t]{0.48\textwidth}
        \centering
        \definecolor{mycolor1}{rgb}{0.36471,0.52157,0.76471}%
\definecolor{mycolor2}{rgb}{0.31373,0.71373,0.58431}%
\definecolor{mycolor3}{rgb}{0.97255,0.72941,0.23529}%
\definecolor{mycolor4}{rgb}{0.91373,0.31373,0.24314}%
\definecolor{mycolor5}{rgb}{0.30100,0.74500,0.93300}%

\begin{tikzpicture}

\begin{semilogyaxis}[xlabel=$k$, ylabel=$\kappa_2$, 
    width=\linewidth,height=5cm,ytick={1,1e25,1e50,1e75,1e100,1e125}]
\addplot [color=mycolor1, line width=0.5pt, mark=*, mark size=1.5pt]
table[x index=0,y index= 1,col sep=comma] {img/condition_numbers_unscaled.txt};
\end{semilogyaxis}

\end{tikzpicture}
        \caption{Condition number {$\kappa_2(G_k)$} of the unscaled $2\times 2$ systems {$G_k$ from \eqref{eq:TwoByTwoConst}--\eqref{eq:TwoByTwoSin}} {for the model with radii} $R_1 = {\SI{44.3e-3}{m}}$, $R_2={\SI{44.7e-3}{m}}$.}
        \label{fig:kappa_unscaled}
    \end{subfigure}
    \hfill
    \begin{subfigure}[t]{0.48\textwidth}
        \centering
        \definecolor{mycolor1}{rgb}{0.36471,0.52157,0.76471}%
\definecolor{mycolor2}{rgb}{0.31373,0.71373,0.58431}%
\definecolor{mycolor3}{rgb}{0.97255,0.72941,0.23529}%
\definecolor{mycolor4}{rgb}{0.91373,0.31373,0.24314}%
\definecolor{mycolor5}{rgb}{0.30100,0.74500,0.93300}%

\begin{tikzpicture}

\begin{semilogyaxis}[xlabel=$k$, ylabel=$\kappa_2$, 
    width=\linewidth,height=5cm]
\addplot [color=mycolor1, line width=0.5pt, mark=*, mark size=1.5pt]
table[x index=0,y index= 2,col sep=comma] {img/condition_numbers.txt};
\end{semilogyaxis}

\end{tikzpicture}
        \caption{Condition number {$\kappa_2(\tilde{G_k})$} of the scaled $2\times 2$ systems {$\tilde{G}_k$ from \eqref{eq:gustav}} {for the model with radii} $R_1=1$, $R_2=1.009$.}
        \label{fig:kappa_scaled}
    \end{subfigure}
    \caption{Comparison of condition numbers for unscaled and scaled systems.}
    \label{fig:condition_number_comparison}
\end{figure}
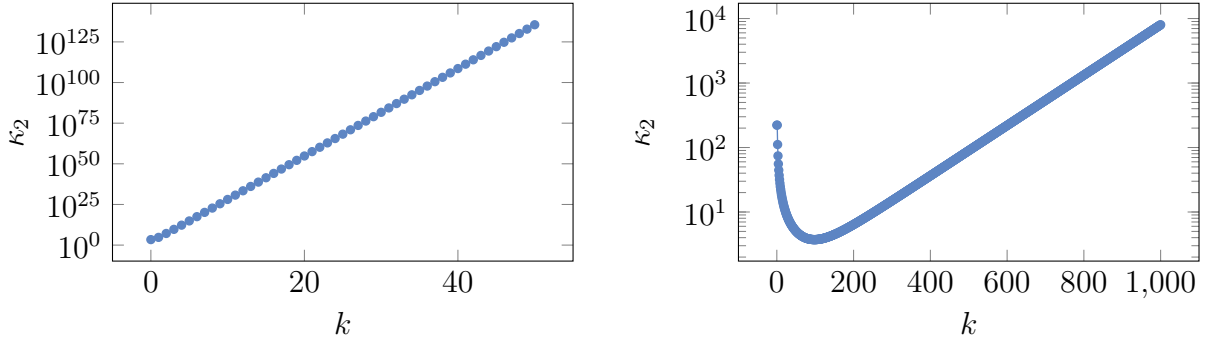

For the computation of the coefficients of the analytical solution \eqref{eq:ana_sol}, one has to solve $2\times 2$ systems 
\eqref{eq:TwoByTwoConst}--\eqref{eq:TwoByTwoSin}, 
where the occurring matrices only depend on $R_1, R_2$ and the mode number $k$. In the later used machine model the radii are $R_1 = {\SI{44.3e-3}{m}}$ and $R_2 = {\SI{44.7e-3}{m}}$. For these values the condition number $\kappa_2$ of the resulting $2\times 2$ systems grows prohibitively as depicted in Figure~\ref{fig:kappa_unscaled}. Therefore, we choose to evaluate the air-gap element on a scaled domain with $R_1 = 1$ and use the scaling laws for the magnetic fields to meet the coupling conditions \cite{Stipetic_2015aa}. This leads to the scaled $2\times2$ systems 
\begin{align}
    \tilde{G}_0{\coloneq}\begin{pmatrix}
        1 & 0 \\
        1 & \ln\left( \tilde{R} \right)
    \end{pmatrix} 
    \text{~~and~~}
    \tilde{G}_k{\coloneq}\begin{pmatrix}
        1 & 1 \\
        \tilde{R}^k & \tilde{R}^{-k}
    \end{pmatrix} \label{eq:gustav}
\end{align}
with $\tilde{R} = \frac{R_2}{R_1}$. As {shown} in Figure \ref{fig:kappa_scaled}, the condition number remains reasonable, even for higher mode numbers. 

\section{Space-Time Formulation}\label{sec:time}
In addition to standard time-stepping schemes, we solve \eqref{eq:dae_system} using {a variational approach in time}. For this purpose let $I=(0,T)$ be the time interval of interest, which is divided into $N_t$ (not necessarily equidistant) intervals, yielding a one dimensional mesh $\mathcal{T}_h$ on which we define piecewise linear and continuous ansatz functions $\lbrace \varphi_i(t) \rbrace_{i=0}^{N_t}$ and piecewise constant test functions $\lbrace \psi_j(t) \rbrace_{j=1}^{N_t}$. The solution in space-time is then given by the ansatz
\begin{align}
    \underline{U}(t) = \sum\limits_{i=0}^{N_t} \underline{U}_i  \varphi_i(t).
\end{align}
Inserting this into \eqref{eq:dae_system} and testing with $\psi_j$ yields
\begin{align}
    \sum\limits_{i=0}^{N_t} \left[ \int_0^T \left( \tilde{M}(t) \dot{\varphi}(t) + \tilde{K}(t) \varphi(t) \right) \psi_j(t)\, dt \right] \underline{U}_i \nonumber \\
    = \int_0^T \tilde{\underline{J}}(t) \psi_j(t) \, dt
\end{align}
This is a discontinous Galerkin formulation in time, introduced e.g. in \cite{Neumuller_2016aa}. If we use midpoint rule for quadrature and denote the mesh width of the $i$-th element by $h_i$ and its midpoint by $t_i$, this leads to the block equation system
\begin{align}
    \begin{pmatrix}
        A_1^- & A_1^+ &       &         &    \\
             & A_2^- & A_2^+  &         &    \\
              &       & \ddots & \ddots  &    \\
              &       &        &A_{N_t}^- &A_{N_t}^+ 
    \end{pmatrix}
    \begin{pmatrix}
        \underline{U}_0 \\
        \underline{U}_1 \\
        \vdots \\
        \underline{U}_{N_t}
    \end{pmatrix}
    = \begin{pmatrix}
        \underline{\tilde{J}}_1 \\
        \underline{\tilde{J}}_2 \\
        \vdots \\
        \underline{\tilde{J}}_{N_t}
    \end{pmatrix} \label{eq:sysInhom}
\end{align}
with 
\begin{align}
    A_i^- &= -\tilde{M}(t_i) + \frac{h_i}{2} \tilde{K}(t_i), \\
    A_i^+ &=  \tilde{M}(t_i) + \frac{h_i}{2} \tilde{K}(t_i), \\
    J_i   &= h_i \tilde{J}(t_i).
\end{align}
Note, that the block matrix \eqref{eq:sysInhom} is rectangular having more columns than rows. Given an initial condition $\underline{U}_0$ the homogenized system is given by
\begin{align}
        \begin{pmatrix}
        A_1^+ &       &         &    \\
        A_2^- & A_2^+  &         &    \\
              & \ddots & \ddots  &    \\
              &        &A_{N_t}^- &A_{N_t}^+ 
    \end{pmatrix}
    \begin{pmatrix}
        \underline{U}_1 \\
        \underline{U}_2 \\
        \vdots \\
        \underline{U}_{N_t}
    \end{pmatrix}
    = \begin{pmatrix}
        \underline{\tilde{J}}_1 - A_1^- \underline{U}_0 \\
        \underline{\tilde{J}}_2 \\
        \vdots \\
        \underline{\tilde{J}}_{N_t}
    \end{pmatrix} \label{eq:sys_initialCond}.
\end{align}
Note that this system is sub-diagonal and hence can be solved via forward substitution. This would correspond to a time stepping scheme. However, the interpretation as a space-time system opens possibility for parallel solvers in the temporal dimension, e.g. via multigrid \cite{Neumuller_2016aa}.
\section{Torque Calculation}\label{sec:torque}
{In electric machines, one of the key performance indicators is the electromagnetic torque, which describes the rotational force generated by the machine. The torque acting on a volume can be determined by integrating the Maxwell stress tensor over a surface that encloses this volume \cite{Salon_1995aa}. For electric machines in a 2D context, the torque is given by}
\begin{align}
T = \frac{r^2 L}{\mu_0} \int_0^{2\pi} B_r B_\varphi \, d \varphi 
\end{align}
with the machine length $L$, the vacuum permeability $\mu_0$ {and the radius of the integration surface/line $r$}. 
{In conventional machine simulations using low-order finite elements, torque evaluation often suffers from inaccuracies. These issues primarily arise from the numerical differentiation of $A_z$ to obtain the required components of $B$ which amplifies numerical errors in the vector potential. The angular component $B_\varphi$  is especially prone to inaccuracies because its computation involves the normal derivative of $A_z$ across element boundaries. To reduce these errors, various techniques have been proposed, such as the eggshell method or Arkkio's method, which use a volume integral instead of a surface integral (or a surface integral in place of a line integral in 2D)  \cite{Arkkio_1987aa,Henrotte_2004aa,Henrotte_2010aa}. Using the air-gap method resolves these accuracy challenges, as the solution in the air-gap is known in closed form, facilitating the derivation of the magnetic vector potential and eliminating errors from numerical differentiation.} 
The $B$ field {can then be} calculated from $A_z$ via
\begin{align}
    B = \begin{pmatrix}
        0 & 1 \\
        -1 & 0
    \end{pmatrix} \nabla A_z
\end{align}
where the matrix entries are with respect to Cartesian coordinates. The gradient in polar coordinates is given~by
\begin{align}
    \nabla A_z = \partial_r A_z \underline{e}_r + \frac{1}{r} \partial_\varphi A_z \underline{e}_\varphi
\end{align}
and straightforward calculations yield the components of the $B$ field in polar coordinates
\begin{align}
    B_r = \frac{1}{r} \partial_\varphi A_z, \quad B_\varphi = -\partial_r A_z.
\end{align}
For any radius $r$ inside the air-gap~\cite{Abdel-Razek_1981aa}, this yields
\begin{align}
    B_r =& \frac{1}{r} \sum\limits_{k=1}^K k \cos(k\varphi) \left( \beta_k r^k + \beta_k' r^{-k} \right) - \frac{1}{r} \sum\limits_{k=1}^K k \sin(k\varphi) \left( \alpha_k r^k + \alpha_k' r^{-k} \right)  \\
    B_\varphi =& -\frac{\A_0'}{r} - \frac{1}{r} \sum\limits_{k=1}^K k \cos(k\varphi) \left( \alpha_k r^k - \alpha_k' r^{-k} \right) - \frac{1}{r}\sum\limits_{k=1}^K k \sin(k\varphi) \left( \beta_k r^k - \beta_k' r^{-k} \right).
\end{align}
Using orthogonality of the occurring functions and {exploiting}
\begin{align}
    \int_0^{2\pi} \sin^2(k \varphi) \, d \varphi = \int_0^{2\pi} \cos^2(k \varphi) \, d \varphi = \pi
\end{align}
yields
\begin{align}
    T &= \frac{\pi L}{\mu_0} \sum\limits_{k=1}^K k^2 \big[
    -\left(\beta_k r^k + \beta_k' r^{-k} \right)
    \left(\alpha_k r^k - \alpha_k' r^{-k} \right) \nonumber \\  
    & \quad \quad \quad \quad \quad \quad + \left(\alpha_k r^k + \alpha_k' r^{-k} \right)\left(\beta_k r^k - \beta_k' r^{-k} \right) \big] \nonumber \\
    &= \frac{2 \pi L}{\mu_0} \sum\limits_{k=1}^K k^2 \left[ \alpha_k'\beta_k - \alpha_k \beta_k' \right].
\end{align}
\section{Machine Model}\label{sec:machine}
\begin{figure}
    \centering
    \scalebox{0.8}{\input{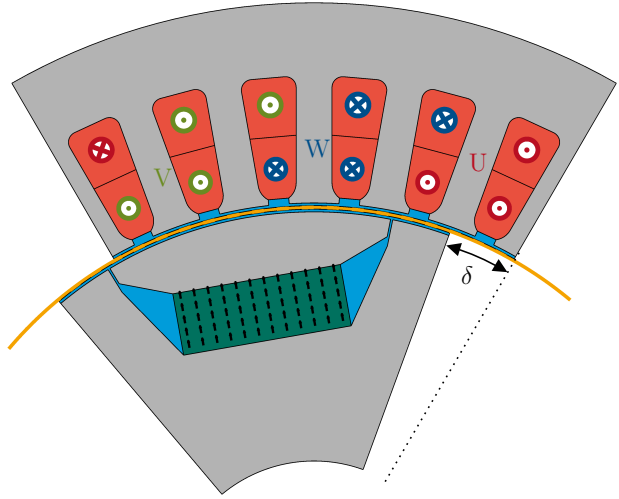}}    
    \caption{Geometry model for the PSMS. The air-gap element connects the current-driven stator to the rotor, where the rotation is prescribed.}
    \label{fig:GeometryMotor}
\end{figure}

To validate the described approach, we consider the permanent magnet synchronous motor (PMSM) given in Figure \ref{fig:GeometryMotor}. 
The model consists of a three-phase driven stator with homogenized copper coils in red, the linear iron yoke in gray, and air shown in blue. The rotor comprises a permanent magnet in green with additional air slits to compensate for nonlinear saturation effects. {This justifies to consider only linear material behavior in the following.} The air-gap element for the coupling is shown as orange band that connects rotor and stator. Due to symmetry, only one sixth of the motor is simulated with appropriate boundary conditions. Eddy currents are assumed to be fully prevalent in the magnet, and in a reduced way in the iron due to lamination which is modeled using the formulation from \cite{Gyselinck_1999aa}. 

The current density is given by the three-phase current $J_z = J_0 \sum_{k=1}^3 \sin\left(\omega t+2/3\,\pi k \right)$ with the pole-pair number $p=3$ and $\omega = 2\pi\cdot\SI{50}{\hertz}$. For synchronous operation, the rotation angle $\delta=\frac{\omega t}{p}$ is prescribed by a constant rotation velocity. 

The rotor and stator domains are discretized using quadratic B-Splines with 6044 and 5196 degrees of freedom, respectively. The analytical solution of the air-gap element is expressed using the first 35 nonzero terms in (\ref{eq:ana_sol}), resulting in 140 unknowns for the coupling. Since only one sixth of the motor is simulated and the solution is antiperiodic, this corresponds to the modes $k\in\{3, 9, 15,..., 207\}$.

\section{Results}\label{sec:results}
\begin{figure}
    \centering
%
%
\definecolor{mycolor1}{rgb}{0.36471,0.52157,0.76471}%
\definecolor{mycolor2}{rgb}{0.31373,0.71373,0.58431}%
\definecolor{mycolor3}{rgb}{0.97255,0.72941,0.23529}%
\definecolor{mycolor4}{rgb}{0.91373,0.31373,0.24314}%
\definecolor{mycolor5}{rgb}{0.30100,0.74500,0.93300}%
\begin{tikzpicture}

\begin{axis}[%
    width=\linewidth,
    height=5.5cm,
    xmin=0,
    xmax=20,
    xlabel style={font=\color{white!15!black}},
    xlabel={$\text{Rotation angle (}^\circ\text{)}$},
    ymin=1.2,
    ymax=1.6,
    ylabel style={yshift=-0.3cm,font=\color{white!15!black}},
    ylabel={Torque (Nm)},
    xtick={0,2,...,20},
    legend pos=south east,
    legend cell align={left},
    legend style={font=\small},
    axis background/.style={fill=white},
    axis x line*=bottom,
    axis y line*=left,
    xmajorgrids,
    ymajorgrids
] 
\addplot [color=black, line width=0.5pt, mark=*, mark size=1.5pt]
table[x index=0,y index= 1,col sep=comma] {img/Torque_ag.txt};
\addlegendentry{Static Reference}

\addplot [color=mycolor1, line width=0.5pt, mark=square*, mark size=1.5pt]
table[x index=0,y index= 2,col sep=comma] {img/Torque_ag.txt};
\addlegendentry{Implicit Euler}

\addplot [color=mycolor4, line width=0.5pt, mark=diamond*, mark size=1.5pt]
table[x index=0,y index= 1,col sep=comma] {img/st_torque.txt};
\addlegendentry{Space-Time}

\end{axis}
\end{tikzpicture}%
    \caption{Torque evaluation for different scenarios. The torque curves are compared for the static case, one time stepping scheme (implicit Euler) and the space-time solution.}
    \label{fig:torque-results}
\end{figure}
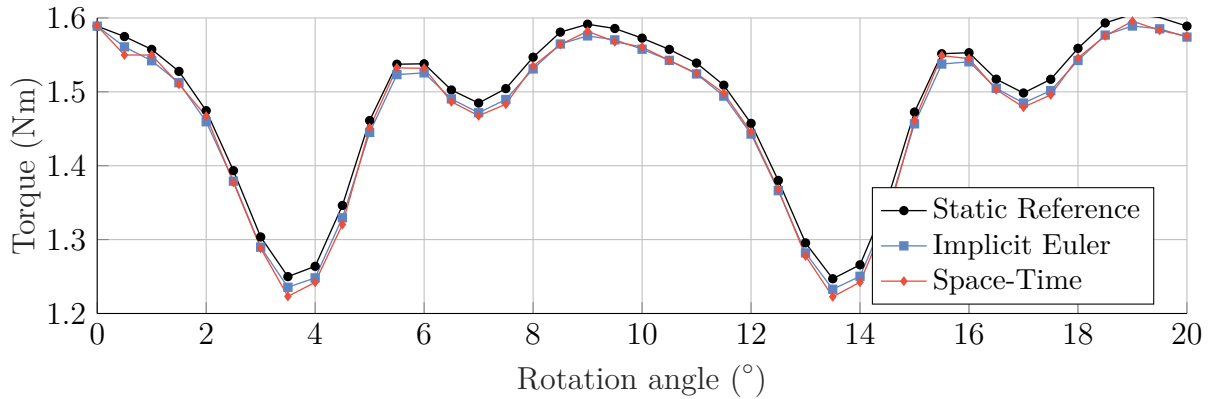

We perform transient simulations on the presented machine for the first $20^\circ$ rotational degrees. After that, the torque will be periodic in the stationary case. The time interval is discretized with $N_t=40$ elements. This corresponds to 40 time-steps of a classical time-stepping scheme. At this point, the resulting block matrix system is solved using a direct solver.
The initial conditions for the problem are computed by a static simulation at time zero. The static and transient simulation results are shown in Figure~\ref{fig:torque-results}. 
The torque values from the static solution are given as reference. As expected, the transient torque is slightly reduced because of the energy losses due to eddy-currents in the magnet and iron which are only captured in the transient case.
The solution from the space-time setting is in very good agreement with the solutions obtained from the time stepping scheme.

\section{Conclusions and Outlook}\label{sec:conclusions}
In this contribution we coupled the analytical solution of the air-gap element to numerical solutions discretized with IGA and solved the transient problem with a space-time formulation. The combination of these methods is advantageous for several reasons: First, the boundary integrals are evaluated precisely due to the exact geometry representation. Second, the rotor movement can be incorporated without the need of remeshing or reevaluation of boundary integrals. Third, the discretization of the time-domain can be {carried out} a priori even without known rotor movement. Forth, the space-time setting allows for new possibilities regarding the parallel solving of the equation system, e.g., using Parareal \cite{Schops_2018aa,Bast_2020aa}, 
and multigrid techniques \cite{Bolten_2020aa}.

Future steps include the extension of the simulations to the nonlinear case, i.e., nonlinear material behavior and unknown rotor motion. {Finally, benchmarking must show if the space-time approach has computational benefits.}

\section*{Acknowledgments}
The work is supported by the joint DFG/FWF Collaborative Research Centre CREATOR (CRC – TRR361 / 10.55776/F90) at TU Darmstadt, TU Graz and JKU Linz.

\end{document}